\documentclass{amsart}
\usepackage{amssymb}
\usepackage{amsmath}
\usepackage{a4}

\def\Gr{\operatorname{Gr}}
\begin{document}
\def\id{\operatorname{id}}
\def\gronk{\vphantom{\vrule height 12pt}}
\def\mapright#1{\smash{\mathop{\longrightarrow}\limits\sp{#1}}}
\newtheorem{theorem}{Theorem}[section]
\newtheorem{lemma}[theorem]{Lemma}
\newtheorem{remark}[theorem]{Remark}
\newtheorem{definition}[theorem]{Definition}
\newtheorem{corollary}[theorem]{Corollary}
\newtheorem{example}[theorem]{Example}
\def\qedbox{\hbox{$\rlap{$\sqcap$}\sqcup$}}
\makeatletter
  \renewcommand{\theequation}{%
   \thesection.\alph{equation}}
  \@addtoreset{equation}{section}
 \makeatother
\def\BB{\mathcal{B}}
\def\MM{{\mathfrak{M}}}
\title[Puffini--Videv models and manifolds]
{Puffini--Videv models and manifolds}
\author{P. Gilkey, E. Puffini, and V. Videv${}^1$}
\begin{address}{PG: Mathematics Department, University of Oregon, Eugene, OR 97403, USA. Email:
\it gilkey@uoregon.edu.}\end{address}
\begin{address}{EP: K. I. T., Malvinas. Email: \it ekaterinapuffin@yahoo.com.}\end{address}
\begin{address}{VV: Veselin Videv, Mathematics Department, Thracian University, University Campus,
6000 Stara Zagora Bulgaria. Email:\it videv@uni-sz.bg.}\end{address}
\begin{abstract} Let $\mathcal{J}(\pi)$ be the higher order Jacobi operator. We study algebraic curvature
tensors where $\mathcal{J}(\pi)\mathcal{J}(\pi^\perp)=\mathcal{J}(\pi^\perp)\mathcal{J}(\pi)$. In the Riemannian
setting, we give a complete characterization of such tensors; in the pseudo-Riemannian setting, partial results
are available. We present non-trivial geometric examples of Riemannian manifolds
with this property.
\end{abstract}
\keywords{$0$-model, algebraic curvature tensor, Einstein manifold, Higher order Jacobi operator, Jacobi operator, Ricci operator.
\newline 2000 {\it Mathematics Subject Classification.} 53C20\newline ${}^1$Corresponding author}
\maketitle

\section{Introduction} The study of commutativity and of spectral properties for natural operators in differential geometry
has received much attention in recent years. Probably the seminal paper in the subject is due to Osserman \cite{O90} who proposed a
characterization of Riemannian rank $1$-symmetric spaces in terms of the spectrum of the Jacobi operator. There are, however, many
other crucial works which should be cited  -- papers by Ivanova and Stanilov \cite{IS95}, by Stanilov
\cite{S04}, by Stanilov and Videv \cite{SV04}, by Szab\'o \cite{Sz91}, and by Tsankov \cite{SV92} are central. However, as the
literature is a vast one, we must limit ourselves and shall refer to the bibliographies in
\cite{GF02,Gl02} for further information.

We shall work in both the geometric and in the algebraic contexts; the Jacobi operator will form the focus of our study. We begin by
introducing some notational conventions.  We say that $\mathfrak{M}:=(V,\langle\cdot,\cdot\rangle,A)$ is a $0$-model if
$\langle\cdot,\cdot\rangle$ is a non-degenerate inner product of signature $(p,q)$ on a finite dimensional 
vector space $V$ of dimension $m=p+q$ and if $A\in\otimes^4V^*$ is an algebraic curvature tensor, i.e. if $A$ is a
$4$-tensor which satisfies the symmetries of the Riemann curvature tensor:
\begin{eqnarray*}
&&A(v_1,v_2,v_3,v_4)=A(v_3,v_4,v_1,v_2)=-A(v_2,v_1,v_3,v_4),\\
&&A(v_1,v_2,v_3,v_4)+A(v_2,v_3,v_1,v_4)+A(v_3,v_1,v_2,v_4)=0\,.
\end{eqnarray*}
The associated curvature operator $\mathcal{A}$ and Jacobi operator $\mathcal{J}$ are then characterized, respectively,
by the identities:
$$\langle\mathcal{A}(v_1,v_2)v_3,v_4\rangle=A(v_1,v_2,v_3,v_4)\quad\text{and}\quad
  \langle\mathcal{J}(v_1)v_2,v_3)=A(v_2,v_1,v_1,v_3)\,.$$
The Jacobi operator $\mathcal{J}(v)$ is quadratic in $v$. It is convenient to polarize and to set
$$
\mathcal{J}(v_1,v_2):v_3\rightarrow\textstyle\frac12\{\mathcal{A}(v_1,v_2)v_3+\mathcal{A}(v_1,v_3)v_2\}\,.
$$
This operator was first introduced to study the geometry of the Jacobi operator by Videv
\cite{Videv93}. If $\vec v:=(v_1,...,v_k)$ is a basis for a non-degenerate $k$-plane $\pi$, let $\xi_{ij}:=\langle v_i,v_j\rangle$
give the components of the metric restricted to $\pi$ relative to the given basis.
If $\xi^{ij}$ denotes the inverse matrix, then one defines the {\it higher order Jacobi operator} by setting:
\begin{equation}\label{eqn-1.a}
\mathcal{J}(\vec v):=\sum_{i=1}^k\sum_{j=1}^k\xi^{ij}\mathcal{J}(v_i,v_j).
\end{equation}
Let $\pi(\vec v):=\operatorname{Span}\{v_1,...,v_k\}$. Then
$\mathcal{J}(\pi):=\mathcal{J}(\pi(\vec v))$
is independent of the particular basis chosen for $\pi$; this operator was first introduced in this context by Stanilov and Videv
\cite{SV92} in the Riemannian setting and latter extended by Gilkey, Stanilov, and Videv to the pseudo Riemannian
setting \cite{GSV98}. Note that $\rho:=\mathcal{J}(V)$ is the {\it Ricci operator};
$$\rho:y\rightarrow\sum_{1=1}^m\sum_{j=1}^m\xi^{ij}\mathcal{A}(y,e_i)e_j\,.$$
Thus the higher order Jacobi operator can also be thought of as a generalization of the Ricci operator to lower
dimensional subspaces.

Let $\Gr_{r,s}(V,\langle\cdot,\cdot\rangle)$ be the Grassmannian of all non-degenerate linear
subspaces of $V$ which have  signature
$(r,s)$; the pair $(r,s)$ is said to be {\it admissible} if and only if $Gr_{r,s}(V,\langle\cdot,\cdot\rangle)$ is non-empty and does
not consist of a single point or, equivalently, if the inequalities $0\le r\le p$, $0\le s\le q$, and $1\le r+s\le m-1$ are satisfied.
Let
$[A,B]:=AB-BA$ denote the commutator of two linear maps. We shall establish the following result in Section
\ref{sect-2}:

\begin{theorem}\label{thm-1.1}
Let $\MM=(V,\langle\cdot,\cdot\rangle,A)$ be a $0$-model. The following assertions are equivalent; if any is satisfied,
then we shall say that $\MM$ is a {\rm Puffini--Videv} $0$-model.\footnote{This notation was suggested by the first author}
\begin{enumerate}
\item There exists $(r_0,s_0)$ admissible so that 
\newline\qquad$\mathcal{J}(\pi)\mathcal{J}(\pi^\perp)=
 \mathcal{J}(\pi^\perp)\mathcal{J}(\pi)$ for all $\pi\in Gr_{r_0,s_0}(V,\langle\cdot,\cdot\rangle)$.
\item 
$\mathcal{J}(\pi)\mathcal{J}(\pi^\perp)=
 \mathcal{J}(\pi^\perp)\mathcal{J}(\pi)$ for every non-degenerate subspace $\pi$.
\item $[\mathcal{J}(\pi),\rho]=0$ for every non-degenerate subspace $\pi$.
\end{enumerate}
\end{theorem}

We say that $\MM=(V,\langle\cdot,\cdot\rangle,A)$ is {\it
decomposible} if there exists a non-trivial orthogonal decomposition
$V=V_1\oplus V_2$ which decomposes $A=A_1\oplus A_2$; in this setting, we shall write $\MM=\MM_1\oplus\MM_2$ where
$\MM_i:=(V_i,\langle\cdot,\cdot\rangle|_{V_i},A_i)$. One says that $\MM$ is {\it
indecomposible} if
$\MM$ is not decomposible. We say
$\MM$ is {\it Einstein} if the Ricci operator $\rho$ is a scalar multiple of the identity. By Theorem \ref{thm-1.1}, any Einstein
$0$-model is Puffini--Videv. More generally, the direct sum of Einstein Puffini--Videv models is again Puffini--Videv; the converse
holds in the Riemannian setting:

\begin{theorem}\label{thm-1.2}
Let $\MM=(V,\langle\cdot,\cdot\rangle,A)$ be a Riemannian $0$-model. Then $\MM$ is Puffini--Videv if and only
if $\MM=\MM_1\oplus...\oplus\MM_k$ where the $\MM_i$ are Einstein.
\end{theorem}

In the pseudo-Riemannian setting, a somewhat weaker result
can be established. One says that a $0$-model is {\it pseudo-Einstein} either if the Ricci operator $\rho$ has only one real eigenvalue
$\lambda$ or if the Ricci operator $\rho$ has two complex eigenvalues $\lambda_1,\lambda_2$ with
$\bar\lambda_1=\lambda_2$. This does not imply that $\rho$ is diagonalizable in the higher signature setting and hence $\MM$ need not
be Einstein.

\begin{theorem}\label{thm-1.3}
 Let $\MM=(V,\langle\cdot,\cdot\rangle,A)$ be a $0$-model of arbitrary signature. If $\MM$ is Puffini--Videv, then we may decompose
$\MM=\MM_1\oplus...\oplus\MM_k$ as the direct sum of pseudo-Einstein $0$-models $\MM_i$.
\end{theorem}

We restrict to the Riemannian context henceforth. Theorem \ref{thm-1.2} yields the characterization that indecomposible Riemannian
$0$-model is Einstein if and only if it is Puffini--Videv. However, in passing to the geometric situation, things become a bit more
complicated. Let
$\mathcal{M}=(M,g)$ be a Riemannian manifold of dimension $m\ge3$; any Riemann surface is automatically Puffini--Videv by Theorem
\ref{thm-1.1}. Let
$\MM_k(\mathcal{M},P)=(T_PM,\langle\cdot,\cdot\rangle,R_P,\nabla R_P,...,\nabla^kR_P)$ be a purely algebraic object which encodes the
geometric information about $\mathcal{M}$ up to the
$(k+2)^{\operatorname{th}}$ order jets of the metric at a point $P$; $\MM(\mathcal{M},P)=\MM_0(\mathcal{M},P)$. The notion of
indecomposibility of a $k$-model is defined as above. One says that
$\mathcal{M}$ is Puffini--Videv if
$\MM(\mathcal{M},P)$ is Puffini--Videv at each point of the manifold. One says $\mathcal{M}$ is {\it locally reducible} at a point $P\in
M$ if there exists a neighborhood $\mathcal{O}$ of $P$ so that $(\mathcal{O},g_M)=(\mathcal{O}_1\times \mathcal{O}_2,g_1\oplus g_2)$
decomposes as a Cartesian product. We say
$\mathcal{M}$ is {\it locally irreducible} at $P$ if this does not happen. Clearly if $\MM_k(\mathcal{M},P)$ is indecomposible for some
$k$, then $\mathcal{M}$ is locally irreducible at $P$. Let 
$$\tau_{\mathcal{M}}=\sum_{ijkl}g^{jk}g^{il}R_{ijkl}$$ 
be the scalar curvature of $\mathcal{M}$. One says
$\mathcal{M}$ exhibits {\it scalar curvature blowup} if there is a geodesic in $\mathcal{M}$ defined on a finite interval $(0,T)$ so
that
$\lim_{t\rightarrow0}|\tau_{\mathcal{M}}(\gamma(t))|=\infty$. Such a manifold is necessarily geodesically incomplete and can not be
embedded isometrically in a geodesically complete manifold.

One has the following examples as we shall discuss further in Section
\ref{sect-3}:

\begin{theorem}\label{thm-1.4}
Let $\mathcal{N}:=(N,ds^2_{\mathcal{N}})$ be a Riemann surface. Assume $\tau_{\mathcal{N}}(P_0)\ne2$ for some $P_0\in N$. Let
$(t,x_1,x_2)$ be local coordinates on
$M:=(0,\infty)\times N$. Give $M$ the metric $g_{\mathcal{M}}$ with non-zero components:
$$g_{\mathcal{M}}(\partial_t,\partial_t)=1\quad\text{and}\quad
g_{\mathcal{M}}(\partial_{x_i},\partial_{x_j})=t^2g_{\mathcal{N}}(\partial_{x_i},\partial_{x_j})\,.$$
\begin{enumerate}
\item $\mathcal{M}:=(M,g_M)$ is a
$3$-dimensional Riemannian Puffini--Videv manifold. 
\item $\MM_1(\mathcal{M},P)$ is indecomposible and $\mathcal{M}$ is locally irreducible at $P=(t,P_0)$.
\item $\mathcal{M}$ exhibits scalar curvature blowup and is not Einstein.
\end{enumerate}
\end{theorem}

\begin{theorem}\label{thm-1.5} 
Let $(x_1,x_2,x_3,x_4)$ be coordinates on $M:=(0,\infty)\times(0,\infty)\times\mathbb{R}^2$. Let $\beta>0$. Give $M$ the metric
$g_\beta$ whose non-zero components are:
\begin{eqnarray*}
&&g_\beta(\partial_{x_1},\partial_{x_1})=g_\beta(\partial_{x_2},\partial_{x_2})=1,\quad
  g_\beta(\partial_{x_3},\partial_{x_3})=x_1^2,\\
&&g_\beta(\partial_{x_4},\partial_{x_4})=x_1(x_1+\beta x_2)\,.
\end{eqnarray*}
\begin{enumerate}
\item $\mathcal{M}_\beta:=(M,g_\beta)$ is a $4$-dimensional
Riemannian Puffini--Videv manifold. 
\item $\MM_2(\mathcal{M}_\beta,P)$ is indecomposible and $\mathcal{M}$ is locally ireducible for all $P\in M$.
\item $\mathcal{M}_\beta$ is not locally isometric to
$\mathcal{M}_\gamma$ for
$\beta\ne\gamma$.
\item $\mathcal{M}_\beta$ exhibits scalar curvature blowup and is not Einstein.
\end{enumerate}
\end{theorem}

\section{Puffini--Videv $0$-models}\label{sect-2}

We begin our study in the algebraic context by establishing Theorem \ref{thm-1.1}. Let $\rho$ be the Ricci operator. Since
$\rho=\mathcal{J}(\pi)+\mathcal{J}(\pi^\perp)$, we have:
\begin{equation}\label{eqn-2.a}
\begin{array}{l}[\rho,\mathcal{J}(\pi)]=
 \{\mathcal{J}(\pi)+\mathcal{J}(\pi^\perp)\}\mathcal{J}(\pi^\perp)
 -\mathcal{J}(\pi^\perp)\{\mathcal{J}(\pi)+\mathcal{J}(\pi^\perp)\}\\
=\mathcal{J}(\pi)\mathcal{J}(\pi^\perp)+\mathcal{J}(\pi^\perp)\mathcal{J}(\pi^\perp)-\mathcal{J}(\pi^\perp)\mathcal{J}(\pi)
-\mathcal{J}(\pi^\perp)\mathcal{J}(\pi^\perp)\\
=[\mathcal{J}(\pi^\perp),\mathcal{J}(\pi)]\,.
\end{array}\end{equation}
This establishes the equivalence of Assertions (2) and (3) in Theorem \ref{thm-1.1}. It is immediate that Assertion (2) implies
Assertion (1). We complete the proof of Theorem \ref{thm-1.1} by showing Assertion (1) implies Assertion (2). Assume there exists
$(r_0,s_0)$ admissible so that
\begin{equation}\label{eqn-2.b}
\mathcal{J}(\pi)\mathcal{J}(\pi^\perp)=\mathcal{J}(\pi^\perp)\mathcal{J}(\pi)\ \forall\ \pi\in
  Gr_{r_0,s_0}(V,\langle\cdot,\cdot\rangle)\,.
\end{equation}

Let $1\le\kappa:=r_0+s_0<m:=\dim(V)$. Let $\{e_1,...,e_{\kappa},e_{\kappa+1},...,e_m\}$ be an orthonormal basis for $V$ where
$\{e_1,...,e_{\kappa}\}$ spans a non-degenerate plane $\pi$ of signature $(r_0,s_0)$. Let $\varepsilon_i:=\langle e_i,e_i\rangle$. Then
$$
\mathcal{J}(\pi):=\sum_{i=1}^\kappa\varepsilon_i\mathcal{J}(e_i)\,.
$$
We
distinguish two cases. Suppose first that $\varepsilon_1=\varepsilon_{\kappa+1}$. Set
$$e_1(\theta):=\cos(\theta)e_1+\sin(\theta)e_{\kappa+1}\,.$$
Then $\{e_1(\theta),e_2,...,e_{\kappa}\}$ is an orthonormal basis for a non-degenerate plane $\pi(\theta)$
of signature $(r_0,s_0)$. One has
\begin{eqnarray*}
0&=&[\rho,\mathcal{J}(\pi(\theta))-\mathcal{J}(\pi)]=0\\
&=&[\rho,(\cos^2\theta-1)\mathcal{J}(e_1)+2\sin\theta\cos\theta\mathcal{J}(e_1,e_{\kappa+1})+\sin^2\theta\mathcal{J}(e_{\kappa+1})]\,.
\end{eqnarray*}
This identity for all $\theta$ implies 
$$
[\rho,\mathcal{J}(e_1)-\mathcal{J}(e_{\kappa+1})]=0\quad\text{if}\quad\varepsilon_1=\varepsilon_{\kappa+1}\,.
$$
Suppose next that $\varepsilon_1=-\varepsilon_{\kappa+1}$. Set 
$e_1(\theta):=\cosh(\theta)e_1+\sinh(\theta)e_{\kappa+1}$.
A similar computation, after paying attention to the signs involved, yields:
$$0=[\rho,(\cosh^2\theta-1)\mathcal{J}(e_1)-2\sinh\theta\cosh\theta\mathcal{J}(e_1,e_{\kappa+1})+\sinh^2\mathcal{J}(e_{\kappa+1})$$
which yields the identity 
$$0=[\rho,\mathcal{J}(e_1)+\mathcal{J}(e_{\kappa+1})]\,.$$
We combine these two calculations to see that for all $1\le i,j\le m$ we have that
\begin{equation}\label{eqn-2.c}
\varepsilon_i[\rho,\mathcal{J}(e_i)]=\varepsilon_j[\rho,\mathcal{J}(e_j)]\,.
\end{equation}

We use Equation (\ref{eqn-2.c}) to see that
$$0=[\rho,\mathcal{J}(\pi)]=\sum_{i=1}^\kappa\varepsilon_i[\rho,\mathcal{J}(e_i)]=\kappa\varepsilon_1[\rho,\mathcal{J}(e_1)]$$
and thus $[\rho,\mathcal{J}(e_1)]=0$. This shows that $[\rho,\mathcal{J}(v)]=0$ for every unit spacelike vector if $s_0>0$ and for
every unit timelike vector if $r_0>0$. We can rescale to conclude $[\rho,\mathcal{J}(v)]=0$ on a non-empty open subset of $V$ and hence,
as this is a polynomial identity, conclude $[\rho,\mathcal{J}(v)]=0$ for all $v\in V$. It then follows from Equation
(\ref{eqn-1.a}) that
$[\rho,\mathcal{J}(\pi)]=0$ for every non-degenerate $k$-plane $\pi$. This completes the proof of Theorem \ref{thm-1.1}. \hfill\qedbox

\medbreak Let $\mathfrak{M}$ be a Puffini--Videv $0$-model. In the pseudo-Riemannian setting, the Ricci operator need
not be diagonalizable. However, we can take the Jordan decomposition to decompose:
$$V=\oplus_\lambda V_\lambda\quad\text{and}\quad\rho=\oplus_\lambda\rho_\lambda$$
where one restricts to $\lambda$ with non-negative imaginary parts and where $\rho_\lambda$ has only the eigenvalue $\lambda$ on
$V_\lambda$ if
$\lambda$ is real and the eigenvalues $\{\lambda,\bar\lambda\}$ on $V_\lambda$ if $\lambda$ is complex. As $[\mathcal{J}(x),\rho]=0$,
$$\mathcal{J}(x)V_\lambda\subset V_\lambda\ \forall\ x,\lambda\,.$$

Let $x_i\in V_{\lambda_i}$ and let $\xi$ be arbitrary. If $\lambda_1\ne\lambda_4$, then $\mathcal{J}(\xi)x_i\in V_{\lambda_1}$. Since
$V_{\lambda_1}\perp V_{\lambda_4}$, $\langle\mathcal{J}(\xi)x_1,x_4\rangle=0$. Setting $\xi=x_2+\varepsilon x_3$ and letting
$\varepsilon$ vary then yields
$\langle\mathcal{J}(x_2,x_3)x_1,x_4\rangle=0$. Consequently
\begin{equation}\label{eqn-2.d}
A(x_1,x_2,x_3,x_4)=-A(x_1,x_3,x_2,x_4)\quad\text{if}\quad\lambda_1\ne\lambda_4\,.
\end{equation}
Suppose that $\lambda_1\ne\lambda_4$ and that $\lambda_2\ne\lambda_4$. We may then compute
\medbreak\qquad
$\phantom{=-}A(x_1,x_2,x_3,x_4)=-A(x_1,x_3,x_2,x_4)$\hfill (Equation (\ref{eqn-2.d}) as
     $\lambda_1\ne\lambda_4$)\smallbreak\qquad
$=\phantom{-}A(x_3,x_2,x_1,x_4)+A(x_2,x_1,x_3,x_4)$\hfill (the Bianchi identity)\smallbreak\qquad
$=-A(x_2,x_3,x_1,x_4)+A(x_2,x_1,x_3,x_4)$\hfill (curvature symmetries)\smallbreak\qquad
$=\phantom{-}A(x_2,x_1,x_3,x_4)+A(x_2,x_1,x_3,x_4)$\hfill (Equation (\ref{eqn-2.d}) as $\lambda_2\ne\lambda_4$)
  \smallbreak\qquad
$=-2A(x_1,x_2,x_3,x_4)$\hfill(curvature symmetries).
\medbreak\noindent This shows 
\begin{equation}\label{eqn-2.e}
A(x_1,x_2,x_3,x_4)=0\quad\text{if}\quad\lambda_1\ne\lambda_4\text{ and }\lambda_2\ne\lambda_4\,.
\end{equation}

Suppose that $A(x_1,x_2,x_3,x_4)\ne0$ and that $\lambda_2\ne\lambda_4$. Then we may use Equation (\ref{eqn-2.e})
to see that $\lambda_1=\lambda_4$ and $\lambda_2=\lambda_3$. Since $\lambda_2\ne\lambda_4$, we may apply Equation
(\ref{eqn-2.d}) to see
$$A(x_1,x_2,x_3,x_4)=-A(x_2,x_1,x_3,x_4)=A(x_2,x_3,x_1,x_4)\,.$$
This vanishes by Equation (\ref{eqn-2.e}) since $\lambda_2\ne\lambda_4$ and $\lambda_3\ne\lambda_4$ which is a
contradiction.

Consequently $A(x_1,x_2,x_3,x_4)\ne0$ implies $\lambda_1=\lambda_2=\lambda_3=\lambda_4$. Thus we may
decompose
$$A=\oplus_\lambda A_\lambda\quad\text{for}\quad A_\lambda\in\otimes^4V_\lambda^*\,.$$
This completes the proof of Theorem \ref{thm-1.3}.
If $\MM$ is Riemannian, then necessarily $\lambda$ is real. Since $\rho$ is self-adjoint with respect to a positive definite metric,
$\rho$ is diagonalizable. This implies $\rho$ is a scalar multiple of the identity and hence $\MM$ is Einstein. Theorem \ref{thm-1.2}
now follows.\hfill\qedbox

\section{Irreducible Riemannian Puffini--Videv non-Einstein manifolds}\label{sect-3}

In the geometric setting, matters are a bit more complicated. The work of Tsankov \cite{T05} can be used to construct
examples of $3$-dimensional manifolds which are irreducible, which are Puffini--Videv, and which are not
Einstein. This construction has been generalized by Brozos-V\'azquez et al. \cite{BG06}. Let $\mathcal{M}$ be as in Theorem
\ref{thm-1.4}. Let $(x_1,x_2)$ be local coordinates on $N$ so the metrics take the form
$$
  ds^2_{\mathcal{N}}=e^{2\alpha}(dx_1^2+dx_2^2)\quad\text{and}\quad
  ds^2_{\mathcal{M}}=dt^2+t^2e^{2\alpha}(dx_1^2+dx_2^2)\,.
$$
The curvature tensor of $\mathcal{N}$ is then given by:
$$
R_{\mathcal{N}}(\partial_{x_1},\partial_{x_2},\partial_{x_2},\partial_{x_1})
=e^{2\alpha}\left(\frac{\partial^2\alpha}{\partial x_1^2}+\frac{\partial^2\alpha}{\partial x_2^2}\right)$$
and work of \cite{BG06}  shows the only non-zero component of the curvature tensor is:
$$R_{\mathcal{M}}(\partial_{x_1},\partial_{x_2},\partial_{x_2},\partial_{x_1})
  =-t^{-2}e^{2\alpha}\left(\frac{\partial^2\alpha}{\partial x_1^2}+\frac{\partial^2\alpha}{\partial x_2^2}+e^{2\alpha}\right)$$
Let $V_1:=\operatorname{Span}\{\partial_{x_1},\partial_{x_2}\}$ and $V_2:=\operatorname{Span}\{\partial_t\}$. We
then have an orthogonal direct sum decomposition $TM=V_1\oplus V_2$. Furthermore, $A=A_1\oplus A_2$ where $A_2$ is trivial
and $A_1$ is Einstein. Thus $\mathcal{M}$ is Puffini--Videv. 

The scalar curvatures on $\mathcal{M}$ and on $\mathcal{N}$ are related by the identity:
$$\tau_{\mathcal{M}}=t^{-2}(\tau_{\mathcal{N}}-2)\,.$$
The curves $t\rightarrow (t,P_0)$ are unit speed geodesics and clearly $\tau_{\mathcal{M}}$ blows up as $t\rightarrow0$. Furthermore
$\mathcal{M}$ is not Einstein since $\tau_{\mathcal{N}}-2$ does not vanish identically. 
The $1$-model is indecomposable since
$\operatorname{Range}\{\mathcal{R}_{\mathcal{M}}\}=\operatorname{Span}\{\partial_1,\partial_2\}$ and since
$\tau_{\mathcal{M}}$ exhibits non-trivial dependence on $t$. This completes the proof of Theorem \ref{thm-1.4}.\hfill\qedbox

\medbreak The proof of Theorem \ref{thm-1.5} is similar. Again, the only non-vanishing curvature is
$R(\partial_{x_3},\partial_{x_4},\partial_{x_4},\partial_{x_3})$. If we set $V_1:=\operatorname{Span}\{\partial_{x_1},\partial_{x_2}\}$
and $V_2:=\operatorname{Span}\{\partial_{x_3},\partial_{x_4}\}$, then $TM=V_1\oplus V_2$ is an orthogonal direct sum. We set $A_1=0$
and let $A_2=R|_{V_2}$. We then have $R=A_1\oplus A_0$ so $\mathcal{M}$ is Puffini-Videv.
The scalar curvature of $\mathcal{M}$ is:
$$\tau_{\mathcal{M}}=x_1^{-1}(x_1+\beta x_2)^{-1}\,.$$
 From this it follows that $\tau_{\mathcal{M}}$
exhibits blowup along the geodesic
$t\rightarrow(t,0,0,0)$ as $t\downarrow0$. Let
$\mathcal{E}:=\operatorname{Range}(\mathcal{R})=\operatorname{Span}\{\partial_{x_1},\partial_{x_2}\}$. Then
$\mathcal{E}^\perp=\operatorname{Span}\{\partial_{x_1},\partial_{x_2}\}$. Let $\Psi:=-\ln|\tau_{\mathcal{M}}|$. The Hessian of $\Psi$
restricted to
$\mathcal{E}^\perp$ takes the form:
\begin{eqnarray*}
H&=&\left(\begin{array}{rr}
-x_1^{-2}+(x_1+\beta x_2)^{-2}&\beta(x_1+\beta x_2)^{-2}\\
\beta(x_1+\beta x_2)^{-2}&\beta^2(x_1+\beta x_2)^{-2}\end{array}\right)
\end{eqnarray*}
It now follows that
$$\det(H|_{\mathcal{E}})=\textstyle\frac14\beta\tau_{\mathcal{M}}^2\,.$$
This shows that $\beta$ is an isometry invariant; in particular $\mathcal{M}_\beta$ is not locally isometric to $\mathcal{M}_\gamma$
if $\beta\ne\gamma$. Since $\det(H|_{\mathcal{E}})$ has rank $2$, it follows easily that the $2$-model is indecomposable. We refer to
\cite{BG06} for further information concerning
the geometry of these manifolds which were first discovered in a different context.\hfill\qedbox

\section*{Acknowledgments}
The research of P. Gilkey was
partially supported by the Max Planck Institute for the Mathematical Sciences (Leipzig, Germany). The research of E. Puffini was
supported by a grant from the K. I. T. It is a pleasure to acknowledge helpful conversations with C. Dunn and with Z. Zhelev
concerning these matters.

\end{document}